\documentclass{article}
\usepackage[utf8]{inputenc}
\usepackage{amsmath, amssymb}
\usepackage{amsthm}
\usepackage{graphicx}

\theoremstyle{definition}
\newtheorem{dfn}{Definition}

\theoremstyle{plain}
\newtheorem{thm}{Theorem}

\theoremstyle{plain}

\theoremstyle{plain}
\newtheorem{lem}{Lemma}

\theoremstyle{definition}

\theoremstyle{definition}
\newtheorem{assump}{Assumption}

\usepackage{biblatex}
\bibliography{bibliography}

\usepackage{apptools}
\AtAppendix{\counterwithin{lem}{subsection}}
\AtAppendix{\counterwithin{thm}{subsection}}

\title{Spectral Bayesian Network Theory}
\author{Luke Duttweiler, Sally W. Thurston, Anthony Almudevar}
\date{\today}

\begin{document}

\maketitle

\begin{abstract}
    A Bayesian Network (BN) is a probabilistic model that represents a set of variables using a directed acyclic graph (DAG). Current algorithms for learning BN structures from data focus on estimating the edges of a specific DAG, and often lead to many `likely' network structures. In this paper, we lay the groundwork for an approach that focuses on learning global properties of the DAG rather than exact edges. This is done by defining the \textit{structural hypergraph} of a BN, which is shown to be related to the inverse-covariance matrix of the network. Spectral bounds are derived for the normalized inverse-covariance matrix, which are shown to be closely related to the maximum indegree of the associated BN. 
\end{abstract}

\textbf{Keywords:} Weighted hypergraph, Bayesian Network, Hypergraph Laplacian, Eigenvalue bound, Directed acyclic graph, Linear structural equation model

\textbf{MSC:} 05C50, 62H22

\section{Introduction}

Bayesian Networks (BNs) are probabilistic models used to model complex systems by representing the relationships among a set of variables with a directed acyclic graph (DAG). BNs are applied in various studies including research on environmental management \cite{barton2012bayesian}, predicting forest fires \cite{sevinc2020bayesian}, and gene regulatory networks \cite{dondelinger2013non}. One of the most important areas of research in Bayesian Networks is structure discovery. This involves determining which edges connect which vertices in the DAG underlying the BN. For this problem there are many existing approaches including the popular PC-stable algorithm \cite{colombo2014order}, Grow-Shrink \cite{margaritis2003learning}, and various Bayesian MCMC techniques which originated with the order-MCMC algorithm \cite{madigan1995bayesian}. 

Given a Bayesian Network with $p$ nodes, structure discovery fundamentally centers on determining which nodes are \textbf{not} connected directly in the DAG. Under common assumptions this is the same as determining which elements of the inverse covariance matrix of the BN are equal to zero. This problem involves the estimation of the $\frac{p(p-1)}{2}$ lower triangular elements of the inverse covariance matrix. As $p$ grows the order of this estimation problem grows quickly.

In order to mitigate this issue, most structure discovery algorithms require an assumption limiting the complexity of the underlying DAG. This most frequently takes the form of an assumption about the maximum number of parents a node can have \cite{almudevar2010hypothesis}. Even when such assumptions are met, if $p$ is large the structure resulting from any particular algorithm is only one among many likely structures, and is unlikely to match the true underlying structure. 

This paper offers a theoretical basis for a method to get around these problems. Our main result, presented as Theorem 1, shows the connection between the inverse covariance matrix of a BN and the Laplcian matrix of a particular weighted hypergraph. The eigenvalues of matrices associated with a hypergraph can be shown to have important meaning in relation to structural properties of the hypergraph (ie. maximum degree, maximum edge-size, etc.). In essence, Theorem 1 allows us to bound structural properties of the BN, using eigenvalues of its inverse-covariance matrix. Since there are only $p$ eigenvalues, any information derived from the eigenvalues can be estimated with a much higher degree of certainty than information derived from the entire inverse covariance matrix.

Theorems 2 and 3 demonstrate more practical uses of Theorem 1, connecting the eigenvalues of the normalized inverse-covariance matrix to the indegree of the associated BN. Both of these theorems suggest direct statistical applications, the details of which will be published in a future paper. 

Section 2 gives a review of the background material needed for the paper, both in hypergraphs and in Bayesian Networks. Section 3 contains Theorem 1, and discusses the relationship between the inverse covariance matrix and a novel construction called the structural hypergraph. In Section 4 we discuss the relationship between eigenvalues of the inverse covariance matrix and the complexity of a BN. Theorems 2 and 3 are presented in Section 4, although their proofs are contained in the Appendix. Section 5 contains a discussion and description of future work. 

\section{Background}

In this section we survey some background material needed for our main results, and define notation that we will use throughout this paper. 

\subsection{Hypergraphs}
 
 We begin by reviewing hypergraphs and their associated matrices. Our definitions for hypergraphs and their matrices are largely based on the definitions provided in Mulas and Reff \cite{mulas2020spectra}, although they are somewhat changed so that we can work with the spectra of weighted hypergraphs in a matrix setting. Galuppi et al. \cite{galuppi2021spectral} approach this topic through tensors, which we do not explore here.
 
\begin{dfn}
A \textit{hypergraph} is defined as the triple $G = (V, E, \mathcal{I})$ where 
\begin{itemize}
    \item $V = \{v_1, \dots, v_{p_V}\}$ is the vertex set,
    \item $E = \{e_1, \dots, e_{p_E}\}$ is the edge set,
    \item $\mathcal{I} \subseteq V\times E$ is the set of incidences.
\end{itemize}

Each edge in $E$ is a subset of $V$. If $v_i$ and $e_k$ are connected within the graph then we say they are \textit{incident}, we have $(v_i, e_k) \in \mathcal{I}$, and $v_i \in e_k$. If $v_i$ and $v_j$ are both incident with $e_k$, and $i \neq j$ then we say $v_i$ and $v_j$ are \textit{adjacent}. In this paper we do not allow a vertex to be adjacent with itself.

A hypergraph is called \textit{incidence-simple} if each vertex is incident with any edge at most once. For the remainder of this paper, all hypergraphs are assumed to be incidence-simple.
\end{dfn}

\begin{dfn}
The \textit{degree} of a vertex, denoted $d_i = deg(v_i)$ is the number of incidences containing $v_i$. The \textit{size} of an edge, denoted $size(e_k)$, is the number of incidences containing $e_k$. We denote the maximum degree and maximum edge size by 

\begin{align*}
    &\Delta = \max_{v\in V} deg(v) \\
    &\nabla = \max_{e\in E} size(e).
\end{align*}
\end{dfn}

\begin{figure}[t]
\center{\includegraphics[width = 90mm]{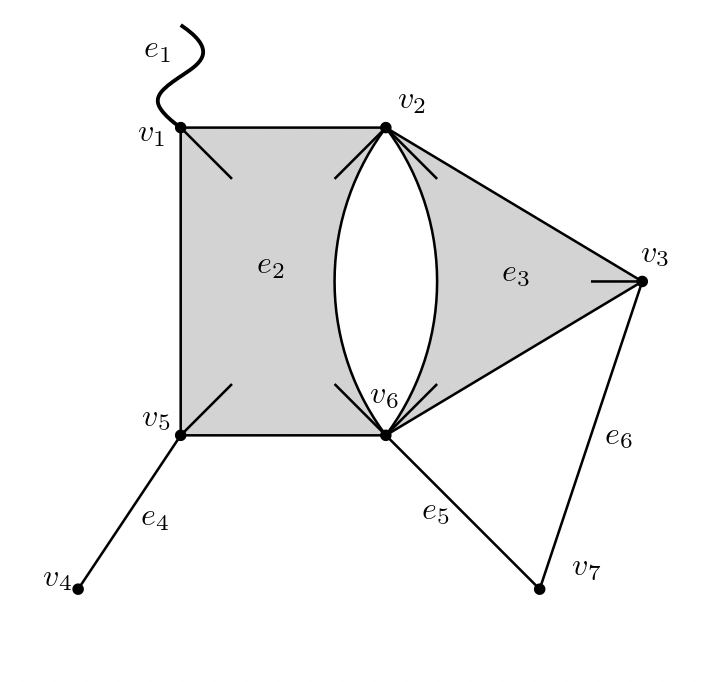}}
\caption{A hypergraph with labeled vertices and edges. Notice that $e_2$ connects 4 vertices and $e_3$ connects 3, while $e_1$ is only connected to one vertex (represented by the edge `dangling'). Additionally, notice that $v_2$ and $v_6$ are adjacent through both $e_2$ and $e_3$.}
\label{Fig1}
\end{figure}

A visual representation of a hypergraph can be seen in Figure \ref{Fig1}.

\begin{dfn}
A \textit{weighted hypergraph} is defined as the quadruple $G = (V, E, \mathcal{I}, \omega),$ where $\omega:V\times E \rightarrow \mathbb{R}$, specifies a \textit{weight} for each vertex-edge incidence for which 

\[
\omega(v_i, e_k) \neq 0 \iff (v_i, e_k) \in \mathcal{I}.
\]
\end{dfn}
\begin{dfn}
The \textit{magnitude} of a vertex $v_i$ and the \textit{effect} of an edge $e_k$ are defined as 

\begin{align*}
    m_i = &mag(v_i) = \sum_{e\in E} \omega(v_i, e)^2 \\
    &eff(e_k) = \sum_{v\in V} \omega(v, e_k)^2.
\end{align*}
\end{dfn}

Notice that magnitude and effect are very similar in concept to degree and size and even exactly the same if the range of $\omega$ is restricted to $\{-1,0,1\}$ as in oriented hypergraphs.

\begin{figure}[t]
    \centering
    \includegraphics[width = 125mm]{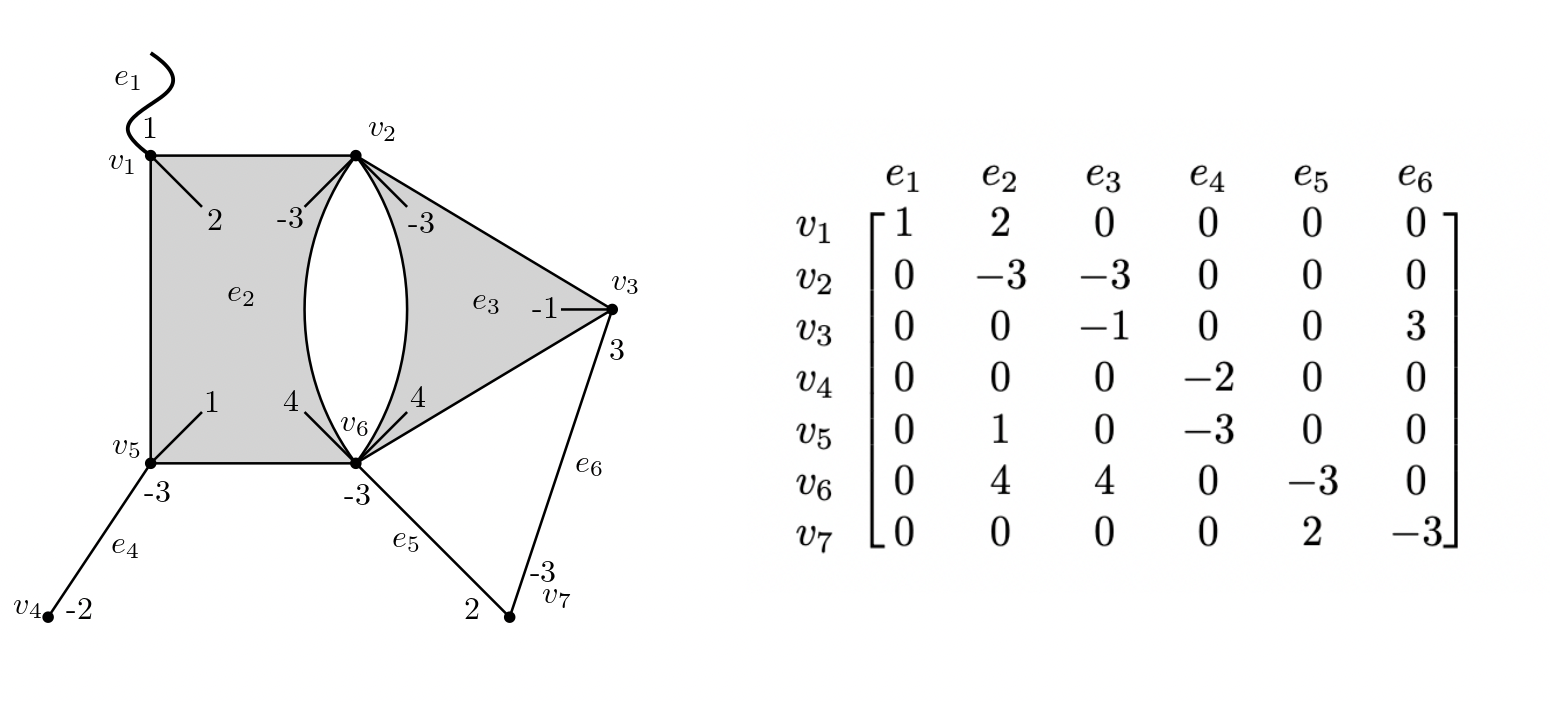}
    \caption{A weighted hypergraph, along with the associated incidence matrix, $H$. Note that while the degree of $v_2$ is 2, the magnitude of $v_2$ is $18.$}
    \label{Fig2}
\end{figure}

\begin{dfn}\label{def: adjWeightFunc}
The \textit{adjacency weight function} is a function $\zeta:E\times V\times V \rightarrow \mathbb{R},$ defined by

\[
\zeta_{e_k}(v_i, v_j) = -\omega(v_i, e_k)\omega(v_j, e_k),
\]

where we set $\zeta_{e_k}(v_i, v_i) = 0$ for all $i, k$.
\end{dfn}
 
 \subsection{Hypergraph Matrices}
The matrix definitions we will use for weighted hypergraphs are again based on the definitions from Mulas (2020), with important differences \cite{mulas2020spectra}.

\begin{dfn}
Let $G = (V, E, \mathcal{I}, \omega)$ be a weighted hypergraph with $|V| = p_V$ and $|E| = p_E$, where $|S|$ denotes the cardinality of a set $S$.

We define the \textit{magnitude matrix} of $G$ by

\[
M(G)^{p_V\times p_V} = diag(m_1, \dots, m_{p_V}).
\]

The \textit{incidence matrix} of $G$ is defined as $H(G)^{p_V \times p_E}$ where

\[
H(G)_{ij} = \omega(v_i, e_j).
\]

The \textit{adjacency matrix} of $G$ is defined as $A(G)^{p_V \times p_V}$ where 

\[
A(G)_{ij} = \sum_{e\in E} \zeta_e(v_i, v_j).
\]

Note that this gives $A(G)_{ij} = 0$ for all $i$.

The \textit{normalized Adjacency matrix} of $G$ is defined as 

\[
C(G)^{p_V\times p_V} = M^{-1/2}(G)A(G)M^{-1/2}(G).
\]

The \textit{Kirchoff Laplacian matrix} of $G$ is defined as

\[
K(G)^{p_V\times p_V} = M(G) - A(G) = H(G)H(G)^T.
\]

The \textit{normalized Laplacian matrix} of $G$ is defined as 

\[
L(G)^{p_V\times p_V} = M^{-1/2}(G)K(G)M^{-1/2}(G) = I_{p_V} - C(G).
\]
\end{dfn}

We label the eigenvalues of a given $p\times p$ symmetric matrix $B$ by 

\[
\lambda_1(B) \geq \dots \geq \lambda_p(B).
\]

All eigenvalues we will examine in this paper are eigenvalues of symmetric matrices, and thus are real numbers. 

Additionally we say that the eigenvalues of a symmetric matrix $B$ are \textit{symmetric about a number} $a \in \mathbb{R}$ if, when $\lambda$ is an eigenvalue of $B$, then so is $2a - \lambda$.
 
\subsection{Linear Bayesian Networks}

We review the definition of a directed acyclic graph and then present a specialized version of a Bayesian Network, similar to the linear SEM found in Loh and Buhlmann (2014).\cite{loh2014high}

We first discuss definitions for direction in a graph. 

\begin{dfn}
Let $G = (V, E)$ be a graph (not a hypergraph, so all edges have size 2). We say an edge $e\in E$ is \textit{directed} if $e = (v_i, v_j) \neq (v_j, v_i)$. Directed edges can be drawn with an arrow showing the direction (ie. $v_i \rightarrow v_j$).

If a sequence of distinct (with the possible exception of the first and last) vertices $(v^*_0, \dots, v^*_s) \subseteq V$ exists such that

\begin{itemize}
    \item $v^*_0 = v_i$
    \item $v^*_s = v_j$
    \item $(v^*_{l-1}, v^*_l) \in E$ for all $l \in 1,\dots,s,$
\end{itemize}

then we say that there is a \textit{path} in $G$ from $v_i$ to $v_j$, with length $s$. If, for all $l \in 1,\dots s,$ $(v^*_{l-1}, v^*_l)$ is a directed edge, then we say that the path is a \textit{directed path}.

We call $G$ a \textit{tree} if there are no paths (directed or undirected) in $G$ that begin and end at the same vertex and $G$ is connected. If this condition holds true but $G$ is not connected, we call $G$ a \textit{sub-tree}.
\end{dfn}

\begin{dfn}
A \textit{directed acyclic graph} (DAG) is defined as a graph $G = (V, E),$ where, 

\begin{itemize}
    \item all $e_k\in E$ are directed
    \item there are no directed paths in $G$ that start and end at the same vertex.
\end{itemize}

Also, we define the set of \textit{parents} of $v_i$ to be 

\[
pa_G(v_i) = \{v\in V| \exists e\in E \text{ with } e = (v, v_i)\},
\]

the set of \textit{children} of $v_i$ to be 

\[
ch_G(v_i) = \{v\in V|\exists e\in E\text{ with } e = (v_i, v)\}.
\]
\end{dfn}

\begin{dfn}
A \textit{weighted directed acyclic graph} (weighted DAG) is defined as the triple $G = (V,E,\xi)$ where $(V,E)$ forms a DAG and $\xi:V\times V \rightarrow \mathbb{R}$ is a function for which

\[
\xi(v_i, v_j) \neq 0 \iff (v_i, v_j) \in E.
\]

If $|V| = p_V$ then we can define a weighted adjacency matrix of a weighted DAG as a $p_V\times p_V$ matrix $A$ such that 

\[
(A)_{ij} = \xi(v_i, v_j).
\]
\end{dfn}

\begin{dfn}
Let $G = (V, E, \xi)$ be a weighted DAG with weighted adjacency matrix $A$, $|V| = p_V$ and $X = (X_v), v\in V$ be a $p_V \times 1$ random vector indexed by $V$. Then, we say that $X$ is a \textit{linear Bayesian Network} (or \textit{linear structural equation model} as in \cite{loh2014high}) with respect to $G$ if 

\[
X = A^TX + \epsilon
\]

where $\epsilon$ is a $p_V \times 1$ vector of independently distributed error terms with $E[\epsilon_i] = 0$ and $Var(\epsilon_i) = \sigma^2_i > 0.$ Additionally for simplicity, and in order to match common statistical notation, we denote the weights by 

\[
(A)_{ij} = \xi(v_i, v_j) = \beta_{ij}.
\]

Notice that if the vertices are ordered such that for all $v_i \in V,$ $pa(v_i) \subseteq \{v_1, \dots, v_{i-1}\},$ then

\[
p(X_i|X_1, \dots, X_{i-1}) = p\big(X_i| \{X_j|v_j \in pa_G(v_i)\}\big).
\]

Thus, a linear Bayesian Network is a specific case of a Bayesian Network as defined by Pearl \cite{pearl2009causality}.

In general we will use $\Sigma$ to denote the covariance matrix of $X$, and 

\[
\Omega := diag(\Sigma^{-1})^{-\frac12}\Sigma^{-1}diag(\Sigma^{-1})^{-\frac12}
\]

to denote the normalized inverse covariance matrix of $X$.

\end{dfn}

\begin{figure}
    \centering
    \includegraphics[width = 125mm]{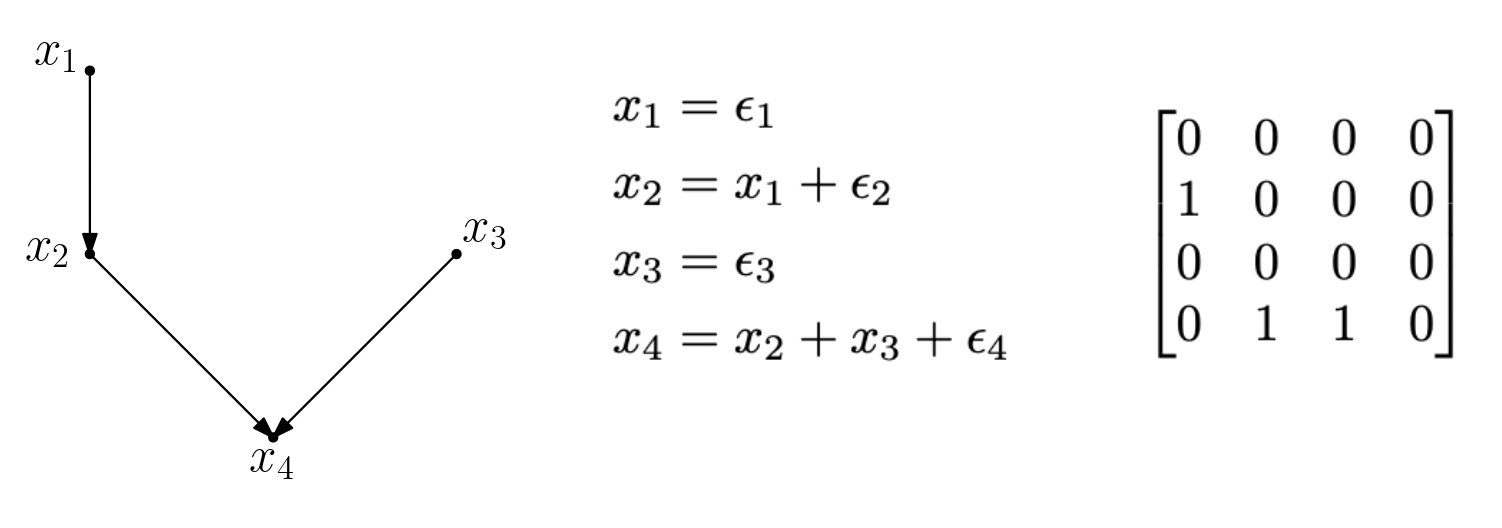}
    \caption{A linear Bayesian network with 4 variables, represented as a DAG, a structural equation model, and the transpose of the adjacency matrix.}
    \label{FigDAG}
\end{figure}

\begin{dfn}
Let $X$ be a linear BN with respect to a weighted DAG $G = (V, E, \xi)$, and let $v_i \in V.$ Following Pearl (1988), we define the \textit{Markov boundary} of $v_i$ to be the union of the parents of $v_i$, the children of $v_i$, and the other parents of the children of $v_i$\cite{pearl1988probabilistic}. That is, 

\[
mb_G(v_i) = pa_G(v_i) \cup ch_G(v_i) \cup \bigg\{v \in pa_G(v_k)|v_k \in ch_G(v_i)\bigg\}\setminus\{v_i\}.
\]

We define the \textit{moral graph} of $X$, denoted $G_M$, as an undirected, unweighted graph (not hypergraph) with adjacency matrix $A_M$ for which

\[
(A_M)_{ij} = 1 \iff v_i \in mb_G(v_j).
\]

This definition follows Frydenberg (1990)\cite{frydenberg1990chain}.
\end{dfn}

\section{The Inverse Covariance Matrix and Hypergraphs}

In this section we define a novel construction on a linear Bayesian Network, which we call the \textit{structural hypergraph}. This weighted hypergraph is then shown to be related to the inverse covariance matrix (frequently used in Bayesian Network structure discovery algorithms).\cite{friedman2008sparse} 

\subsection{Structural Hypergraphs}

\begin{dfn}\label{structuralDef}
Let $X$ be a linear Bayesian Network following a weighted DAG $G$ with weighted adjacency matrix $A$. Also let $|V_G| = p_V$ and $|E_G| = p_E.$ We define the \textit{structural hypergraph of X} to be the weighted hypergraph $G_{ST} = (V_{ST}, E_{ST}, \mathcal{I}_{ST}, \omega_{ST})$ with 

\begin{itemize}
    \item $V_{ST} = V_G,$
    \item $E_{ST} = \{e_1, \dots, e_{p_V}\}$ where $e_k = \{v_k\} \cup pa_G(v_k),$
    \item $\mathcal{I}_{ST}\subseteq V_{ST}\times E_{ST},$ such that  $(v_i, e_k) \in \mathcal{I}_{ST} \iff v_i \in e_k$, 
    \item $\omega_{ST}(v_i, e_k) = \begin{cases}
    -\frac{\beta_{ik}}{\sigma_k} &\text{ if }i\neq k \\
    \frac{1}{\sigma_k} &\text{ if } i=k.
    \end{cases}$
\end{itemize}

Critically, $|E_{ST}| = |V_G| = p_V$.

Also, notice that this gives $(v_i, e_k) \in \mathcal{I}_{ST} \iff \omega_{ST}(v_i, e_k) \neq 0.$
\end{dfn}

With the structural hypergraph defined we are now able to present our primary result. Theorem \ref{invcovStructure} shows the connection between the structural hypergraph and the inverse covariance of its linear Bayesian Network, allowing us to access spectral hypergraph theory in the study of linear Bayesian Networks.

\begin{figure}[t]
    \centering
    \includegraphics[width = 125mm]{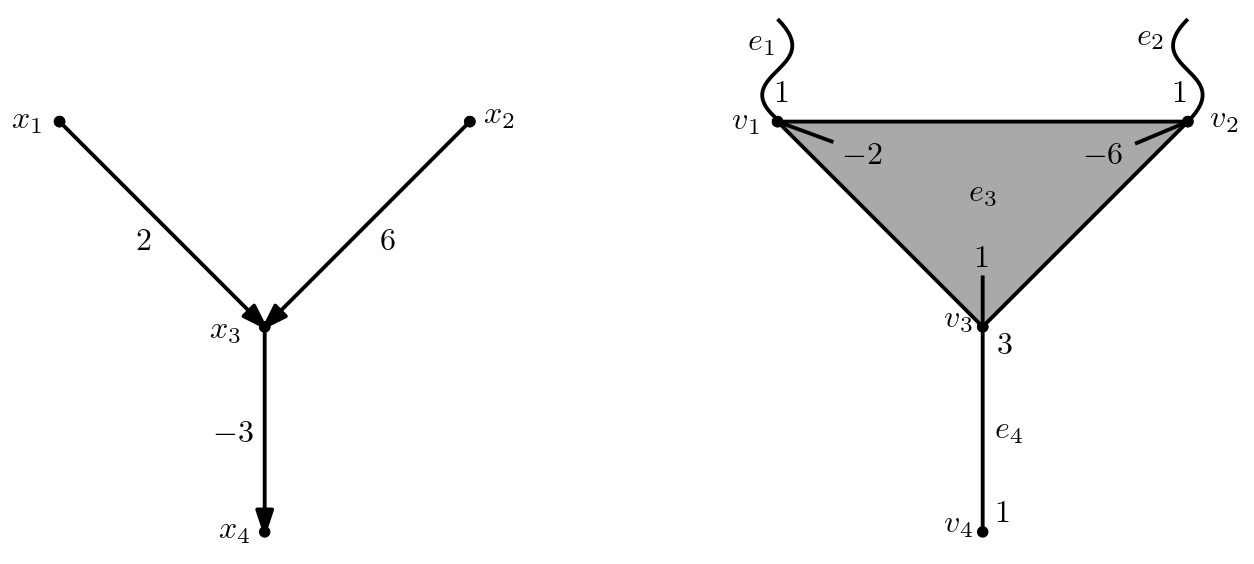}
    \caption{A simple linear Bayesian Network with its structural hypergraph.}
    \label{FigSH}
\end{figure}

\begin{thm}\label{invcovStructure}
Let $X$ be a linear Bayesian Network following a DAG $G$ with weighted adjacency matrix $A$, and let $G_{ST}$ be the structural hypergraph of $X$. Then, if $\Sigma$ is the covariance matrix of $X$, 

\[
\Sigma^{-1} = K(G_{ST}).
\]

Additionally, if the normalized inverse covariance matrix of $X$ is $\Omega$, then

\[
\Omega = L(G_{ST}).
\]
\end{thm}

\noindent\textbf{Proof:}

First, since we have

\[
X = A^TX + \epsilon
\]

we will denote $\Sigma_\epsilon = diag(\sigma^2_1, \dots, \sigma^2_{p_V})$ to be the covariance matrix of $\epsilon$. 

Now, notice that because $A$ is nilpotent, then we must have that $I-A$ is non-singular. Therefore we can say,

\[
X = (I-A^T)^{-1}\epsilon,
\]

which gives 

\[
\Sigma = (I-A^T)^{-1}\Sigma_\epsilon(I-A)^{-1}.
\]

From this we easily see

\[
\Sigma^{-1} = (I-A)\Sigma_\epsilon^{-1}(I-A^T) = \Big((I-A)\Sigma_\epsilon^{-1/2}\Big)\Big((I-A)\Sigma_\epsilon^{-1/2}\Big)^T
\]

Now, for simplicity denote $B = \Big((I-A)\Sigma_\epsilon^{-1/2}\Big)$ and notice that 

\[
(B)_{ik} = \begin{cases}
-\frac{\beta_{ik}}{\sigma_k} &\text{ if }i\neq k \\
\frac{1}{\sigma_k} &\text{ if }i=k.
\end{cases}
\]

Therefore $B = H(G_{ST})$ by definition, and we have

\[
\Sigma^{-1} = BB^T = H(G_{ST})H(G_{ST})^T = K(G_{ST}).
\]

Of course this also easily gives by definition that 

\[
\Omega = L(G_{ST}). \blacksquare
\]

\begin{figure}[t]
    \centering
    \includegraphics[width = 125mm]{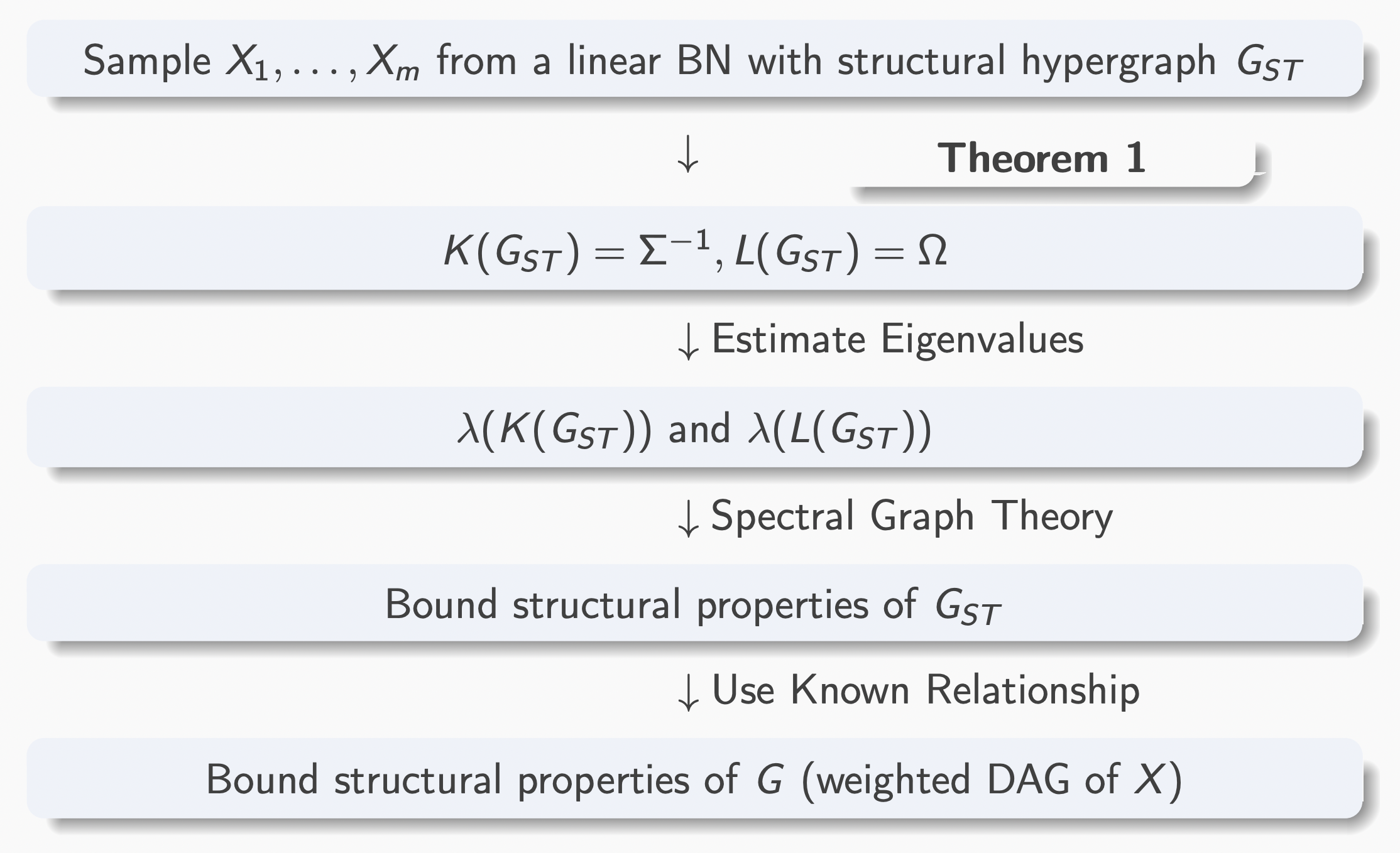}
    \caption{A graphic demonstrating the use of Theorem \ref{invcovStructure}.}
    \label{Fig4}
\end{figure}

Theorem \ref{invcovStructure} points out that the inverse covariance matrix of a linear BN is the Laplacian matrix of a particular weighted hypergraph with structural properties related to the linear BN. Therefore, if we estimate the eigenvalues of the inverse covariance matrix, we can use bounds on the spectra of the structural hypergraph to learn about the structural properties of the respective linear Bayesian Network. Figure \ref{Fig4} gives a visual representation of this idea.

In the remainder of the paper we demonstrate a use for the result of Theorem \ref{invcovStructure}. By leveraging relationships between the eigenvalues of $\Omega = L(G_{ST})$ and the structural properties of $G_{ST}$, we develop theoretical statements about how the true eigenvalues of $\Omega$ relate to the complexity of the respective linear Bayesian Network.

\section{How Eigenvalues Relate to Bayesian Network Complexity}

Under the assumption that the maximum number of parents for any node $K$ is 1 (in which case the moral graph $G_M$ is a tree or sub-tree), Chow and Liu (1968) gives an algorithm for learning a Bayesian Network from data that runs in polynomial time\cite{chow1968approximating}. In more recent years it has become apparent that assuming a maximum on number of parents is critical as allowing $K \geq 3$ gives an NP-Hard problem, and even $K = 2$ can be difficult\cite{chickering1996learning}\cite{chickering2004large}. However, without expert knowledge of the topic at hand, it can sometimes be difficult to justify any assumption about $K$ in a scientifically meaningful way. 

In the following section we provide theoretical results which rely on Theorem \ref{invcovStructure} about the relationship between the eigenvalues of the normalized inverse covariance, $\Omega$, and the maximum number of parents for any node, $K$. In future work we will be developing hypothesis tests based on these relationships. 

\subsection{The Largest Eigenvalue}

Our first result in this section pertains to the maximum eigenvalue of the normalized inverse covariance matrix. The proof of Theorem \ref{treeEigenvalue} (and associated Lemmas) is provided in Appendix A.

\begin{thm}\label{treeEigenvalue}
    Let $X$ be a linear Bayesian Network with moral graph $G_M$. Then if $G_M$ is a tree (or sub-graph of a tree), and $\Omega$ is the normalized precision matrix of $X$ we must have
    
    \[
    \lambda_1\big(\Omega\big) \leq 2.
    \]
\end{thm}

This result suggests a simple hypothesis test to determine whether or not the assumption that $G_M$ is a tree is reasonable given data. A hypothesis test based on Theorem \ref{treeEigenvalue} has the extraordinary property of only requiring the estimation of one value (the largest eigenvalue of $\Omega$) regardless of the number of nodes in the linear BN!

It is important to note a severe limitation of Theorem \ref{treeEigenvalue}. Simply knowing that $\lambda_1(\Omega) \leq 2$ does not guarantee that the moral graph is a tree. This problem is something we seek to overcome with a somewhat more complicated result in the following section.  

\subsection{Eigenvalue Symmetry}

Following two key assumptions, we present Theorem 3 in this section, which related the additive symmetry of the eigenvalues of $\Omega$ to the maximum indegree of a BN. The necessary proofs for this section may be found in Appendix B.

\begin{assump}\label{assump: 1}
Let $X$ be a linear Bayesian Network with structural hypergraph $G_{ST} = (V_{ST}, E_{ST}, \mathcal{I}_{ST}, \omega_{ST}).$ Then, for all $i \neq j$ we have

\[
\sum_{e\in E_{ST}}\omega_{ST}(v_i, e)\omega_{ST}(v_j, e) = 0 \implies \omega_{ST}(v_i, e_k)\omega_{ST}(v_j, e_k) = 0
\]

for all $e_k \in E$.
\end{assump}

\begin{assump}\label{assump: 2}
Let $X$ be a linear Bayesian Network with structural hypergraph $G_{ST}$ and let $s$ be an odd, positive integer. Then, $\text{tr}\big(C(G_{ST})^s\big) = 0$ implies that all diagonal elements of $C(G_{ST})^s$ equal 0.
\end{assump}

Assumption \ref{assump: 1} constitutes a typical probabilistic faithfulness assumption common in the study of Bayesian Networks (see \cite{loh2014high}, Assumption 1). On the other hand, Assumption \ref{assump: 2} is, to the authors' knowledge, a new assumption for linear BNs necessary for our results here. It is crucial to note that if the weights of a given linear Bayesian Network are independently sampled and continuous random variables, then Assumptions 1 and 2 are violated with probability 0. We state this result as a Lemma and provide the proof in Appendix B.

\begin{lem}\label{lemma: assumptions}
Let $X$ be a linear Bayesian Network following a weighted DAG $G$ with weighted adjacency matrix $A$. Let the non-zero weights $(A)_{ij} = \beta_{ij}$ be independently sampled continuous random variables. Then

\begin{itemize}
    \item Assumption \ref{assump: 1} is satisfied except on a set of probability 0 and,
    \item Assumption \ref{assump: 2} is satisfied except on a set of probability 0.
\end{itemize}
\end{lem}

Now, with Assumptions 1 and 2 in place, we are able to give a powerful result. 

\begin{thm}\label{Theorem3}
Let $X$ be a linear Bayesian Network following a DAG $G$ with moral graph $G_M$ and structural hypergraph $G_{ST}.$ Let $\Omega$ be the normalized inverse covariance matrix of $X$. Then the following statements are true:
\begin{itemize}
    \item[(a)] If $G_M$ is a tree, the eigenvalues of $\Omega$ are additively symmetric about 1. 
    \item[(b)] Under Assumptions 1 and 2, if the eigenvalues of  $\Omega$ are additively symmetric about 1, then $G_M$ is a tree.
\end{itemize}
\end{thm}

While the result in Theorem \ref{Theorem3} is somewhat more complex than in Theorem 2, a hypothesis test based on Theorem 3 remains a considerable improvement over attempting to determine whether $G_M$ is a tree by trying to learn which elements of $\Sigma^{-1}$ are non-zero. This is because a test involving the eigenvalues of $\Omega$ requires estimating $p$ values, while $\Sigma^{-1}$ has $p(p-1)/2$ distinct values that must be estimated. 

\section{Discussion and Future Work}

The eigenvalues of the normalized inverse covariance matrix can be estimated with much greater precision and ease than current methods which attempt to discover the entire structure of a Bayesian Network. By making the connection between these eigenvalues and the structural hypergraph in Theorem 1, we provide a means to interpret the eigenvalues in terms of structural properties of the Bayesian Network. Additionally, Theorems 2 and 3 demonstrate actual structural knowledge, which may justify the use of considerably more efficient learning algorithms, gained from simply estimating the eigenvalues rather than attempting to learn the entire structure. 

There are, of course, statistical considerations of importance as the eigenvalues derived from a sample normalized inverse covariance matrix will not adhere strictly to the rules described above. We are currently in the process of developing theory for hypothesis tests based on Theorems 2 and 3, and plan publish the details of these tests soon. 

Additionally, we are hopeful that there are other connections between the eigenvalues of the normalized inverse covariance matrix and the structural properties of a linear Bayesian Network. Spectral Graph Theory has demonstrated many relationships between eigenvalues of a normalized Laplacian matrix and its graph, and it stands to reason that there may be more information to be learned about a BN without having to estimate the entire structure. 

\vspace{.4cm}

\noindent\textbf{Declaration of Competing Interest:}

There is no competing interest.

\vspace{.4cm}

\noindent\textbf{Acknowledgements:}

The authors would like to acknowledge the helpful comments of Dr. Howard Skogman and Dr. Nathan Reff.

Research reported in this publication was supported by the National Institute of Environmental Health Sciences of the National Institutes of Health (NIH) under award number T32ES007271.  The content is solely the responsibility of the authors and does not necessarily represent the official views of the NIH.

\newpage

\appendix

\section*{Appendix A: Theorem 2}

\def\thesection{A}

\subsection{A Weighted Hypergraph Upper Bound}

In this section we derive an upper bound on the spectra of the normalized Laplcian matrix of any weighted hypergraph. This bound was originally demonstrated for oriented hypergraphs in \cite{mulas2020spectra}, and we simply extend it to weighted hypergraphs. For our purposes its primary significance is only as a stepping stone to Theorem \ref{treeEigenvalue}.

We begin by defining the Rayleigh-Ritz quotient and stating the Min-Max theorem, proven in \cite{courant1954methods}.

\begin{lem}\label{MinMax}
    Let $A$ be a Hermitian $n\times n$ matrix and define the Rayleigh-Ritz quotient of $A$ as the function $R_A: \mathbb{R}^{n \times 1} \rightarrow \mathbb{R}$, such that
    
    \[
    R_A(x) = \frac{x^TAx}{x^Tx}.
    \]
    
    Then we have 
    \[
    \lambda_1(A) = \max_{x\in\mathbb{R}^{n \times 1}}\{R_A(x) ; x \neq 0\}
    \]
    
    and 
    
    \[
    \lambda_n(A) = \min_{x\in\mathbb{R}^{n \times 1}}\{R_A(x) ; x \neq 0\}.
    \]
\end{lem}

Lemma \ref{rayleighNorm} gives the exact form of the Rayleigh-Ritz quotient for the normalized Laplacian of a weighted hypergraph.

\begin{lem}\label{rayleighNorm}
    Let $G = (V,E, \mathcal{I}, \omega)$ be a weighted hypergraph with $|V| = p_V$ and $|E| = p_E$. Then 
    
    \[
        R_{L(G)}(x) = \frac{x^TL(G)x}{x^Tx} = \frac{\sum_{e\in E}\bigg(\sum_{i=1}^{p_V}\frac{\omega(v_i, e)}{\sqrt{mag(v_i)}}x_i\bigg)^2}{\sum_{i=1}^{p_V}x_i^2}.
    \]
\end{lem}

\textbf{Proof:}

First, observe that 

\[
    x^TL(G)x = \Big(x^TM^{-1/2}H\Big)\Big(x^TM^{-1/2}H\Big)^T.
\]

Then, since 

\[
    x^TM^{-1/2}H = \bigg[\sum_{i = 1}^{p_V}\frac{\omega(v_i, e_1)}{\sqrt{mag(v_i)}}x_i, \dots, \sum_{i = 1}^{p_V}\frac{\omega(v_i, e_{p_E})}{\sqrt{mag(v_i)}}x_i\bigg],
\]

we must have 

\[
    x^TL(G)x = \sum_{e\in E}\bigg(\sum_{i=1}^{p_V}\frac{\omega(v_i, e)}{\sqrt{mag(v_i)}}x_i\bigg)^2,
\]

giving our result.$\blacksquare$

We now derive an upper bound for the largest eigenvalue of $L(G)$ for any weighted hypergraph $G$.

\begin{lem}\label{maxSmallerThanEdge}
    Let $G = (V,E, \mathcal{I}, \omega)$ be a weighted hypergraph with $|V| = n$. Then, 
    
    \[
    \lambda_1(L(G)) \leq \nabla.
    \]
\end{lem}

\textbf{Proof:}

First, observe that from Lemma \ref{rayleighNorm} we know that for any $x$ such that $||x|| = 1$ we must have 

\[
    R_{L(G)}(x) = \frac{\sum_{e\in E}\bigg(\sum_{i=1}^n\frac{\omega(v_i, e)}{\sqrt{mag(v_i)}}x_i\bigg)^2}{\sum_{i=1}^nx_i^2} = \sum_{e\in E}\bigg(\sum_{i=1}^n\frac{\omega(v_i, e)}{\sqrt{mag(v_i)}}x_i\bigg)^2.
\]

In particular, let $y = (y_1, \dots, y_n)$ be the vector of norm 1 which maximizes $R_{L(G)}(x)$. Then, by Lemma \ref{MinMax} we know that 

\[
\lambda_1\big(L(G)\big) = R_{L(G)}(y).
\]

Now, observe that for all $e\in E$, since in general $a^2 + b^2 \geq 2ab$, we must have 

\begin{align*}
    \Bigg(\sum^n_{i=1}\frac{\omega(v_i, e)}{\sqrt{mag(v_i)}}y_i\Bigg)^2 &= \sum_{v_i\in e}\bigg(\frac{\omega(v_i, e)}{\sqrt{mag(v_i)}}y_i\bigg)^2 + \sum_{v_i, v_j \in e : i\neq j} 2\bigg(\frac{\omega(v_i, e)}{\sqrt{mag(v_i)}}y_i\bigg)\bigg(\frac{\omega(v_j, e)}{\sqrt{mag(v_j)}}y_j\bigg) \\
    &\leq \sum_{v_i\in e}\bigg(\frac{\omega(v_i, e)}{\sqrt{mag(v_i)}}y_i\bigg)^2 + \sum_{v_i, v_j \in e : i\neq j}\bigg(\frac{\omega(v_i, e)}{\sqrt{mag(v_i)}}y_i\bigg)^2 + \bigg(\frac{\omega(v_j, e)}{\sqrt{mag(v_j)}}y_j\bigg)^2 \\.
\end{align*}

Now, in the second sum observe that each term $\Big(\frac{\omega(v_i, e)}{\sqrt{mag(v_i)}}y_i\Big)^2$ will occur for each $v_j \in e$ when $i\neq j.$ Therefore, for each edge $e$,

\begin{align*}
    \Bigg(\sum^n_{i=1}\frac{\omega(v_i, e)}{\sqrt{mag(v_i)}}y_i\Bigg)^2 &\leq \sum_{v_i\in e}\bigg(\frac{\omega(v_i, e)}{\sqrt{mag(v_i)}}y_i\bigg)^2 + \sum_{v_i\in e}\Big(|e| - 1\Big)\bigg(\frac{\omega(v_i, e)}{\sqrt{mag(v_i)}}y_i\bigg)^2 \\ 
    &= |e|\sum_{v_i\in e}\bigg(\frac{\omega(v_i, e)}{\sqrt{mag(v_i)}}y_i\bigg)^2.
\end{align*}

This gives us that 

\begin{align*}
    \lambda_1\big(L(G)\big) &= \sum_{e\in E}\bigg(\sum_{i=1}^n\frac{\omega(v_i, e)}{\sqrt{mag(v_i)}}y_i\bigg)^2 \\
    &\leq \sum_{e\in E} |e|\sum_{v_i\in e}\bigg(\frac{\omega(v_i, e)}{\sqrt{mag(v_i)}}y_i\bigg)^2 \\
    &= \sum_{e\in E} |e|\sum_{i=1}^n\bigg(\frac{\omega(v_i, e)}{\sqrt{mag(v_i)}}y_i\bigg)^2 \\
    &= \sum_{i = 1}^n \sum_{e \in E} |e|\bigg(\frac{\omega(v_i, e)}{\sqrt{mag(v_i)}}y_i\bigg)^2 \\
    &=\sum_{i=1}^n\frac{y_i^2}{mag(v_i)}\sum_{e\in E}|e|\omega(v_i, e)^2 \\ 
    &\leq \nabla\sum_{i=1}^n\frac{y_i^2}{mag(v_i)}\sum_{e\in E}\omega(v_i, e)^2 \\
    &= \nabla\sum_{i=1}^n\frac{y_i^2}{mag(v_i)}mag(v_i) \\
    &= \nabla.
\end{align*}

Therefore, we have that 

\[
\lambda_1\big(L(G)\big) \leq \nabla.\blacksquare
\]

This bound is of great interest as it relates the spectra of a matrix easily derived from a weighted hypergraph to a structural property of that weighted hypergraph which \textbf{does not} depend on the weights.

\subsection{Structural Hypergraph Edge Size in Trees}

\begin{lem}\label{moralTree}
    Let $X$ be a linear Bayesian Network with associated DAG $G$, moral graph $G_M$, and structural hypergraph $G_{ST}$. Then $G_M$ is a tree (or sub-tree) if and only if the maximum edge size in $G_{ST}$ is $\nabla_{ST} \leq 2$.
\end{lem}

\textbf{Proof:} $(\Rightarrow)$ First, assume $G_M$ is a tree (or sub-graph of a tree). Notice that if any vertex in $G$ has more than one parent then the parents are joined by an edge in $G_M.$ This would create an undirected cycle in $G_M$, causing $G_M$ to not be a tree. Therefore,

\[
    \max_{v\in G}\{|pa(v)|\} \leq 1.
\]

Now, observe that in $G_{ST}$ the size of edge $e_i = |pa(v_i)| + 1$. Therefore 

\[
\nabla_{ST} = \max_{v\in G}\{|pa(v)|\} + 1 \leq 2.
\]

$(\Leftarrow)$Now, assume that $\nabla_{ST} \leq 2$. In this case we know that 

\[
\max_{v\in G}\{|pa(v)|\} \leq 1,
\]

and therefore the moral graph $G_M$ does not contain any edges which are not already in $G$. Thus, if $G_M$ contains an undirected cycle, then $G$ must also contain that cycle. However, because $G$ is a DAG, the only way for $G$ to contain an undirected cycle is for at least one vertex in the cycle to have two parents. Therefore, since we know the maximum number of parents is 1 or less, $G_M$ cannot contain an undirected cycle, and is therefore a tree (or sub-graph of a tree).

Thus, $G_M$ is a tree (or sub-graph of a tree) $\iff \nabla_{ST} \leq 2.\blacksquare$

\subsection{Proof of Theorem 2}

Now we are ready to prove Theorem 2. Fortunately, all of the work has already been done, we need to simply use the Lemmas that have been shown. 

\textbf{Proof:}

Let $G_{ST}$ be the structural hypergraph for our linear BN $X$, and let the moral graph $G_M$ be a tree or sub-tree. Then by Lemma \ref{moralTree} we must have that $G_{ST}$ has maximum edge size $\nabla_{ST}\leq 2.$ Therefore, by Lemma \ref{maxSmallerThanEdge} we have $\lambda_1(L(G_{ST})) \leq \nabla_{ST} \leq 2.$ Therefore our result is shown. $\blacksquare$

\section*{Appendix B: Theorem 3}

\def\thesection{B}
\setcounter{subsection}{0}

\subsection{Assumptions 1 and 2}

Assumptions \ref{assump: 1} and \ref{assump: 2} are critical for the proofs in the remainder of Appendix B. However as Lemma \ref{lemma: assumptions} states, if the non-zero weights on the BN are assumed to be drawn independently from continuous distributions then both Assumptions are violated only on sets of probability 0. We provide a proof of Lemma \ref{lemma: assumptions} here. 

\textbf{Proof:}

Let $X$ be a linear BN on $p_V$ variables, following a DAG $G$ with weighted adjacency matrix $A$, and let $G_{ST}$ be the structural hypergraph of $X$. Let the non-zero (as determined by the DAG $G$) weights $(A)_{ij} = \beta_{ij}$ all be independently drawn as continuous random variables. Let $\sigma_1, \dots, \sigma_{p_V}$ be positive random variables sampled independently from the weights. 

We now consider Assumption \ref{assump: 1}. Let $\sum_{e\in E}\omega(v_i, e)\omega(v_j, e) = 0$ for some $i\neq j$. In opposition to Assumption \ref{assump: 1}, assume that there exists some $e_k \in E$ such that $\omega(v_i, e_k)\omega(v_j, e_k) \neq 0.$

Now, WLOG say $j \neq k$ and thus we have by definition

\[
\omega(v_i, e_k)\omega(v_j, e_k) = \begin{cases}
\frac{\beta_{ik}\beta_{jk}}{\sigma_k^2} &\text{ if }i\neq k \\
\frac{\beta_{jk}}{\sigma_k^2} &\text{ if }i = k.
\end{cases}
\]

Thus, since $\sum_{e\in E}\omega(v_i, e)\omega(v_j, e) = 0$ we must be able to write $\beta_{jk}$, a continuous random variable, as a function of random variables of which it is independent. This, of course, can occur only with probability 0. 

We now consider Assumption \ref{assump: 2}. While more significantly more complex, since for all $i \neq j$

\[
    C(G_{ST})_{ij} = -\frac{\sum_{e\in E} \omega_{ST}(v_i, e)\omega_{ST}(v_j, e)}{\sqrt{mag(v_i)mag(v_j)}},
\]

then it is easy to see that if $s$ is a positive, odd integer, then $\text{tr}\Big(C(G_{ST})^s\Big)$ is a sum of terms which are functions of the DAG weights. Therefore, just as for Assumption 1, $\text{tr}\Big(C(G_{ST})^s\Big) = 0$ implies that each diagonal element of $C(G_{ST})^s$ is 0 with probability 1.

Thus, Lemma \ref{lemma: assumptions} is proven. $\blacksquare$

\subsection{The Moral Graph and the Structural Hypergraph}

\begin{lem}\label{moralAdj}
Let $X$ be a linear Bayesian Network following a DAG $G$, with moral graph $G_M$ and structural hypergraph $G_{ST}$. Then, if $A_M$ is the adjacency matrix of $G_M$, where $i \neq j$ we have that $(A_M)_{ij} = 1$ if and only if $\omega_{ST}(v_i, e_k)\omega_{ST}(v_j, e_k) \neq 0$ for at least one $e_k \in E_{ST}.$
\end{lem}

\textbf{Proof:}

Observe that:

\begin{align*}
    (A_M)_{ij} = 1 &\iff v_i \in mb_G(v_j) \\
    &\iff v_i \in pa_G(v_j) \textbf{ or } v_i \in ch_G(v_j) \textbf{ or } \exists  v_l \in V \text{ s.t. } v_i, v_j \in pa(v_l) \\
    &\iff (v_i, e_k) \in \mathcal{I}_{ST} \text{ and } (v_j, e_k) \in \mathcal{I}_{ST},
\end{align*}

where $k = i$, $k = j$, or $k = l$. 

Then, since by definition we have $(v_i, e_k) \in \mathcal{I}_{ST} \iff \omega_{ST}(v_i, e_k) \neq 0$, for the proper value of $k$,

\begin{align*}
    (v_i, e_k) \in \mathcal{I}_{ST} \text{ and } (v_j, e_k) \in \mathcal{I}_{ST} &\iff \omega_{ST}(v_i, e_k)\omega_{ST}(v_j, e_k) \neq 0.\blacksquare
\end{align*}

\begin{lem}\label{moralStr}
Let $X$ be a linear Bayesian Network following a weighted DAG $G$ with moral graph $G_M$ and structural hypergraph $G_{ST}$. Then, if $A_M$  is the adjacency matrix of $G_M$, under Assumption \ref{assump: 1} we have,

\[
(C(G_{ST}))_{ij} = 0 \iff (A_M)_{ij} = 0
\]

\end{lem}

\textbf{Proof:}

Observe that, since

\begin{align*}
    C(G_{ST})_{ij} &= -\frac{\sum_{e\in E} \omega_{ST}(v_i, e)\omega_{ST}(v_j, e)}{\sqrt{mag(v_i)mag(v_j)}}, \\
\end{align*}

using the definition of a structural hypergraph we have,

\begin{align*}
    C(G_{ST})_{ij} = 0 &\iff \sum_{e\in E} \omega_{ST}(v_i, e)\omega_{ST}(v_j, e) = 0. \\
\end{align*}

Thus, by Assumption \ref{assump: 1} we have $\omega_{ST}(v_i, e_k)\omega_{ST}(v_j, e_k) = 0$ for all $e_k \in E$. Since Lemma \ref{moralAdj} gives, 

\[\omega_{ST}(v_i, e_k)\omega_{ST}(v_j, e_k) = 0 \text{ for all }e_k\in E \iff (A_M)_{ij} = 0,\] 

we have our result. $\blacksquare$

\subsection{The Bipartite Moral Graph}

\begin{lem}\label{bipartiteTree}
Let $G = (V,E)$ be a DAG with moral graph $G_M$. Then, $G_M$ is a tree if and only if $G_M$ is bipartite.
\end{lem}

\textbf{Proof:}

Because $G_M$ is an unweighted, undirected graph, it is clear that if $G_M$ is a tree, then $G_M$ is bipartite. We prove the other direction by contradiction.

Let $G_M$ be bipartite and assume that $G_M$ is not a tree. Since $G_M$ is not a tree, there exists an undirected cycle within $G_M$. Because $G_M$ has an undirected cycle there exists at least one vertex $v \in V$ such that $|pa_G(v)| \geq 2.$ Let $v_i, v_j \in pa_G(v).$ Then, by the definition of $G_M$ we know that $(v_i, v), (v_j, v), (v_i, v_j) \in E$, and $G_M$ cannot be bipartite, which is a contradiction. 

Thus

\[
G_M \text{ is a tree,} \iff G_M \text{ is bipartite. }\blacksquare
\]

\subsection{Proof of Theorem 3}

\textbf{Proof:}

First, we prove part (a). Let $G_M$ be a tree. 

By Lemma \ref{bipartiteTree} we then know that $G_M$ is a bipartite graph. Therefore, if $A_M$ is the adjacency matrix of $G_M$ there exists a $r \times q$ matrix $B$ such that 

\[
A_M = \begin{pmatrix}
0_q  & B^T\\ 
B & 0_r \\
\end{pmatrix}
\]

where $0_q, 0_r$ are the zero matrices of size $q\times q$ and $r\times r$ respectively, and $q+r = p_V$. 

Then, by Lemma \ref{moralStr} we know that there exists a $r \times q$ matrix $B_{ST}$ such that 

\[
C(G_{ST}) = \begin{pmatrix}
0_q  & B_{ST}^T\\ 
B_{ST} & 0_r \\
\end{pmatrix}
\]

Therefore 

\[
-C(G_{ST}) = D^{-1}C(G_{ST})D
\]

where $D$ is a diagonal matrix where the diagonal is $q$ 1s followed by $r$ -1s.

Therefore, $C(G_{ST})$ is similar to $-C(G_{ST})$, and by definition, the eigenvalues of $C(G_{ST})$ are symmetric about 0. Then, since by Theorem \ref{invcovStructure}, $I_{p_V} - C(G_{ST}) = \Omega$, we must have that if $\lambda$ is an eigenvalue of $C(G_{ST})$ then $\lambda +1$ is an eigenvalue of $\Omega$. Therefore, by definition, the eigenvalues of $\Omega$ are symmetric about 1. 

\vspace{.2cm}

Now we show part (b). Let the eigenvalues of $\Omega$ be symmetric about 1. Then we know that the eigenvalues of $C(G_{ST})$ must be symmetric about 0. Let $(\lambda_1, \dots, \lambda_n)$ be the spectrum of $C(G_{ST})$. 

Because these eigenvalues are symmetric about 0 we know that 

\begin{align*}
    \lambda_1 + \lambda_2 + &\dots + \lambda_{p_V} = 0 \\
    \lambda_1^3 + \lambda_2^3 + &\dots + \lambda_{p_V}^3 =0 \\
    \lambda_1^5 + \lambda_2^5 + &\dots + \lambda_{p_V}^5 = 0\\
    &\hspace{.2cm}\vdots
\end{align*}

Therefore, we must have 

\begin{align*}
    \text{tr}\Big(C(G_{ST}&)\Big) = 0 \\
    \text{tr}\Big(C(G_{ST}&)^3\Big) = 0 \\
    \text{tr}\Big(C(G_{ST}&)^5\Big) = 0 \\
    &\vdots
\end{align*}

Then, by Assumption \ref{assump: 2}, the diagonal elements of $C(G_{ST}), C(G_{ST})^3, C(G_{ST})^5, \dots$ are all 0. Through Lemma \ref{moralStr}, this implies that, under Assumption \ref{assump: 1}, the diagonal elements of $A_M, A_M^3, A_M^5, \dots$ are all 0. Since $A_M$ is the adjacency matrix of an undirected, unweighted graph $G_M$, this implies that $G_M$ contains no walks of odd length that start and end at the same vertex. Therefore $G_M$ is bipartite, and by Lemma \ref{bipartiteTree}, $G_M$ is a tree.

Therefore, Theorem 3 is proven. $\blacksquare$

\newpage
\printbibliography
\end{document}